\documentclass[letterpaper,conference,10pt]{IWMRRF}
\usepackage{amsmath,amssymb,amsfonts,amscd,amsbsy,amsxtra}
\usepackage{commath} 
\usepackage{graphicx}

\renewcommand{\epsilon}{\varepsilon}
\renewcommand{\d}{\mathrm{d}}

\usepackage{url}       

\usepackage{amsmath}   
\usepackage{array}

\begin{document}

\title{Covariant geometric characterization of slow invariant manifolds: New concepts and viewpoints
}

\author{\authorblockN{
Dirk Lebiedz\authorrefmark{1}}
\authorblockA{\authorrefmark{1}Institute for Numerical Mathematics, Ulm University, Germany}}

\maketitle

\begin{abstract}
We point out a new view on slow invariant manifolds (SIM) in dynamical systems which departs from a purely geometric covariant characterization implying coordinate independency. The fundamental idea is to treat the SIM as a well-defined geometric object in phase space and elucidate characterizing geometric properties that can be evaluated as point-wise analytic criteria. For that purpose, we exploit curvature concepts and formulate our recent variational approach in terms of coordinate-independent Hamiltonian mechanics.

\noindent
Finally, we combine both approaches and conjecture a differential geometric definition of slow invariant manifolds.
For the Davis-Skodje model  the latter can be completely expatiated.
\end{abstract}

\section{Introduction}

\noindent
The analysis of slow invariant manifolds (SIM) for singularly perturbed dynamical systems goes back to Fenichel's geometric singular perturbation theory (GSPT). It formulates by help of perturbation techniques for the critical manifold (singular perturbation parameter $\epsilon=0$) the classical SIM existence theorem involving non-uniqueness similar to center manifold theory.  A coordinate-specific graph representation of the slow manifold (mapping slow to fast variables) allows for a Taylor-series in the singular perturbation parameter $\epsilon$ whose coefficients can be determined successively by matched-asymptotic expansion from the invariance equation.

\noindent
Our aim is to avoid the technical asymptotic local analysis as a first step to address the issue of bundling of trajectories near slow invariant manifolds to be rigorously defined as mathematical objects. Instead, we propose a purely geometric approach that uses global analysis in terms of coordinate-independent differential geometric concepts. 
The first step into that direction goes back to 2004 \cite{Lebiedz2004} when we suggested to formulate a trajectory-based variational principle distinguishing the SIM trajectories - parameterized by fixed reaction progress variables ({\em rpv}) initial values - from other trajectories. The objective functional for this problem formulation has been developed from entropy production \cite{Lebiedz2004} to curvature in time-parametrization \cite{Lebiedz2011a} over the years. 
As in many other computational techniques for point-wise slow manifold approximation, a numerical grid for the {\em rpv} is chosen and the non-reaction progress variables ({\em non-rpv}) are computed on each grid point such that the resulting point in phase space is close to the SIM \cite{Lebiedz2013}.

\begin{figure}[th]
	\centering
	\includegraphics[scale=0.2]{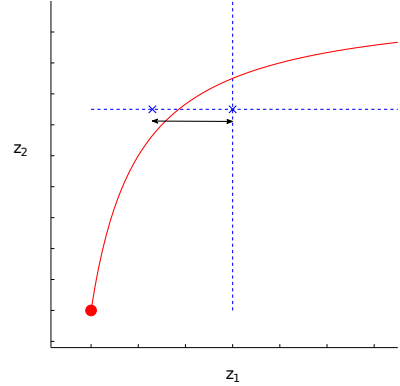}
	\caption{Schematic illustration of the dependency of computational SIM (red curve) approximation (blue crosses) on {\em rpv} choice for the Davis-Skodje example model (\ref{ds-eq}).} \label{fig:covariance}
\end{figure}

\noindent
Most of these numerical SIM approximation methods, in particular the practical ones commonly used in applications, suffer from coordinate-dependency, they yield (at least slightly) different results depending on the chosen coordinate system ({\em rpv}) for the SIM (see Fig. \ref{fig:covariance}). We conjecture that a covariant formulation of SIM criteria, coordinate-indepedency, i.e. invariance of the evaluation result under the {\em rpv} permutation group, should be at the root of a geometric characterization of an object in phase space that is reasonably related to the Fenichel SIM.

\noindent
We follow two roads to address this issue, the first one \cite{Lebiedz2016} is analogy reasoning to analytic mechanics with its coordinate independent general structure allowing to easily figure out conserved properties in the time evolution of a dynamical system. 
The second one  \cite{Heiter2017} is rooted in differential geometry as the theoretical framework for covariant analysis and coordinate-independent formulation of local geometric properties of manifolds and exploits curvature concepts.

\noindent
We develop and illustrate our reasoning on the background of the benchmark problem (Davis-Skodje model)

\begin{eqnarray} \label{ds-eq}
   \dot z_1 & = & -z_1 , \\
          \epsilon \dot z_2 & = & -z + \frac{z_1}{1+z_1} - \frac{\epsilon z_1}{(1+z_1)^2} \nonumber
\end{eqnarray}

\noindent
as an example for a singularly perturbed two time-scale system with analytically known SIM $z_2(z_1)=\frac{1}{1+z_1}$.

\begin{figure}[th]
	\centering
	\includegraphics[scale=0.3]{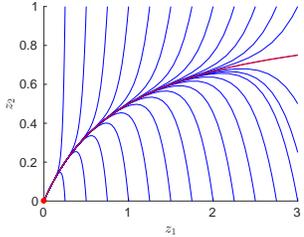}
	\caption{Some solution trajectories (blue) and SIM (red) for the Davis-Skodje model  (\ref{ds-eq}) in phase space (time as curve parametrization), $\epsilon=10^{-2}$.} \label{fig:ds-sym}
\end{figure}

\section{Analytic Mechanics}
Inspired by an obvious local symmetry in the attraction behaviour of solution trajectories to the SIM (see Fig. \ref{fig:ds-sym}, contraction rate ratio of vectors in normal versus tangent bundle of the SIM locally becomes extremal \cite{Adrover2007}) we formulate our variational problem for SIM approximation \cite{Lebiedz2011a} in a Hamiltonian context \cite{Lebiedz2016} in order to be able to identify a conserved property that is related to a symmetry of the Lagrangian by the famous Noether-theorem. 

\noindent
We find an objective functional that allows analytical computation of the exact SIM of the Davis-Skodje model in terms of the solution of a variational boundary value problem \cite{Lebiedz2016}:

\begin{equation}\label{var_exact_ds}
 \min_{z(t)} \int_{t_{0}}^{t_{\rm{f}}}\!k_1\Big\|\frac{\d}{\d t}z(t)\Big\|_2^2-k_2\|z(t)\|_2^2 \; \textrm{d} t
\end{equation}
subject to
\begin{eqnarray}
	 \dot z_1 & = & -z_1 , \nonumber \\
          \epsilon \dot z_2 & = & -z + \frac{z_1}{1+z_1} - \frac{\epsilon z_1}{(1+z_1)^2}, \ z_1(t_f)=z_1^{t_f} \nonumber
\end{eqnarray} 
with constants $k_1,k_2\in\mathbb{R}$. Chosing $k_1=1,k_2=\tfrac{\epsilon^{-1}}{z_1(t)+1}$, the solution of problem (\ref{var_exact_ds}) on arbitrary time horizon $[t_0,t_f]$ yields exactly the SIM $z_2(z_1)=\frac{1}{1+z_1}$.

\noindent
In the Hamiltonian viewpoint model reduction via projection to the SIM can be interpreted as partial integration of the dynamical system flow on the basis of the existence of a first integral (the Hamiltonian is conserved along the flow!).

\section{Differential Geometry}

\noindent
One of the most famous examples of a covariant physical theory is Albert Einstein's general theory of relativity based on a space-time manifold with appropriate metric tensor derived from the gravitational field equation. The appropriate mathematical framework is Riemannian geometry.

\noindent
By analogy reasoning with respect to general relativity we construct the extended phase space by adding to the phase space of a dynamical system the time as a fully equivalent state variable. A particular solution of the differential equation model is then a geometric curve in extended phase space. The SIM is a state-space-time manifold. We conjecture that an appropriate metrization (definition of a Riemannian metric) of the solution manifold in extended phase space will allow a complete geometric characterization of the SIM. The force-curvature relation from general relativity can be transferred to chemical kinetics.
\begin{figure}[th]
	\centering
	\includegraphics[scale=0.25]{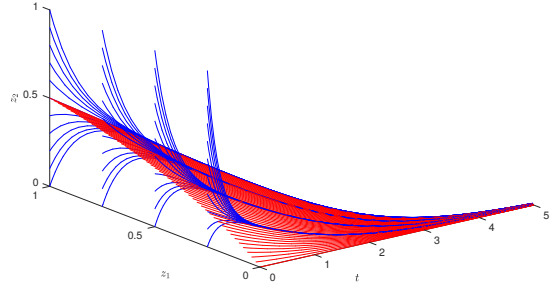}
	\caption{Extended phase space with some solution trajectories (blue) and SIM (red) for the Davis-Skodje model (\ref{ds-eq}), $\epsilon=10^{-2}$.} \label{fig:extendphasespace_DS}
\end{figure}
\subsection{Curvature and Invariance: Necessary Condition for SIM}
\noindent
We demonstrated \cite{Heiter2017} that the SIM invariance property can in general be covariantly formulated as vanishing time-sectional curvatures of the extended-phase-space graph of the SIM $z_2(z_1,t)$ transporting in time along the flow an initial value manifold $z_2(z_1,0)$. Vanishing curvature can be formulated point-wise as a root finding analytic criterion providing a necessary condition for the SIM.

\subsection{Sufficient Condition for SIM: Davis-Skodje Model}
\noindent
Specifically for the Davis-Skodje model  (\ref{ds-eq}), the evaluation of the time-sectional curvature criterion yields a differential equation that can be analytically solved to a parameterized solution of the model equation. In this case the analytic mechanical approach and the covariant differential geometric viewpoint can be combined to formulate a sufficient criterion for the SIM \cite{Heiter2017}.





\end{document}